\documentclass[11pt, reqno]{amsart}
\usepackage{amssymb}
\usepackage{amssymb,amsmath,amscd}
\usepackage[pdfstartview=FitH ]{hyperref} 
\usepackage[hyperpageref]{backref}
\renewcommand{\d}{\delta}

\newcommand{\GL}{\text{\rm GL}}

\newcommand{\A}{\Bbb{A}}

\newcommand{\f}{\frac}

\newcommand{\G}{\Gamma}

\newcommand{\bea}{\begin{eqnarray}}
\newcommand{\eea}{\end{eqnarray}}
\newcommand{\bna}{\begin{eqnarray*}}
\newcommand{\ena}{\end{eqnarray*}}

\renewcommand{\le}{\left}
\newcommand{\ri}{\right}
\newcommand{\s}{\sigma}

\newcommand{\ve}{\varepsilon}

\title{\bf A THEOREM ON ANALYTIC STRONG MULTIPLICITY ONE}

\begin{document}
\openup 1\jot
\maketitle

\centerline{\sc Jianya Liu\footnote{Supported in part by
the 973 Program, and by NSFC Grant \# 10531060.} {\&} Yonghui
Wang\footnote{Supported in part by the 973 Program
\#2007CB807900-2007CB807903, and by NSFC Grant \# 10601034.}}


\bigskip
\bigskip
\begin{quote}
{\footnotesize
\centerline{\scshape Abstract}
Let $K$ be an algebraic number field, and $\pi=\otimes\pi_{v}$ an
irreducible, automorphic, cuspidal representation of $\GL_{m}(\mathbb{A}_{K})$ with analytic
conductor $C(\pi)$. The theorem on analytic strong multiplicity one established in
this note states, essentially, that there exists a positive constant
$c$ depending on $\varepsilon>0, m,$ and $K$ only, such that $\pi$ can be decided completely
by its local components $\pi_{v}$ with norm
$
N(v)<c\cdot C(\pi)^{2m+\varepsilon}.
$

}

\end{quote}

\bigskip

\section{Introduction}
\setcounter{equation}{0}
Let $K$ be an algebraic number field, let
$\pi$ and $\pi'$ be two cuspidal automorphic representations of
$\GL_m(\A_K)$ with restricted tensor product decompositions
$\pi=\otimes \pi_v$ and $\pi'=\otimes\pi'_v$. The strong multiplicity one
theorem states that if $\pi_{v}\cong\pi'_v$ for all but finitely many places $v,$
then $\pi=\pi'.$ The reader is referred to \cite{Brumley2004}
for history and references in this direction.

The analytic version of the above theorem gives, in terms of the
analytic conductor $C(\pi)$ of $\pi$ defined in (\ref{ANAcon}), more precise information on
the number of places needed to decide a cuspidal automorphic representation $\pi$ of
$\GL_m(\A_K)$. Such an analytic result was first established by Moreno
\cite{Moreno1985}, by applying zero-free regions of the Rankin-Selberg $L$-function
of two automorphic representations.
To state the result, let $\mathcal{B}_{m}(Q)$ denote the set of all
cuspidal automorphic representations $\pi$ on $\mathrm{GL}_{m}(\mathbb{A}_{K})$ with
analytic conductors $C(\pi)$ less than a large real number $Q.$ Suppose that
$\pi=\otimes\pi_{v}$ and
$\pi'=\otimes\pi'_{v}$ are in $\mathcal{B}_{m}(Q)$ with
$m\geq 2.$ Then, according to \cite{Moreno1985},
there exist positive constants $c$ and $d$ such that, if
$\pi_{v}\cong\pi'_{v}$ for all finite places $v$ with norm
\bea\label{Moreno}
N(v)\leq
\left\{\begin{array}{lll}
cQ^d, & \text{for } m=2, \\
c\exp(d\log^2Q), & \text{for } m\geq 3,
\end{array}
\right.
\eea
then $\pi=\pi'.$

Using a different method, Brumley \cite{Brumley2004} strengthened Moreno's
result in (\ref{Moreno}) to
\bea\label{Brumley}
N(v)\leq cQ^d, \quad \text{for } m\geq 1,
\eea
where $c$ and $d$ are positive constants depending on
$m$. Brumley's method \cite{Brumley2004} actually proves that
\bea\label{Brumley/C=}
d=\frac{17}{2}m-4+\varepsilon
\eea
is acceptable in (\ref{Brumley}).
By applying a still different method, the second named author
\cite{Wang} showed essentially that, for all $K$, the $d$ in
(\ref{Brumley}) can be reduced to
\bea\label{Wang}
d=4m+\varepsilon.
\eea

The purpose of this note is to show that a suitable modification of the
argument in \cite{Wang} actually gives the acceptable value
\bea\label{Wang}
d=2m+\varepsilon
\eea
in (\ref{Brumley}).
To state the result, let $\mathcal{A}_{M}(Q)$ denote the set of all cuspidal
automorphic representations $\pi$ on $\mathrm{GL}_{m}(\mathbb{A}_{K})$, with
$1\leq m\leq M$,
whose analytic conductors $C(\pi)$ are less than a large real number $Q.$ Thus,
$$
\mathcal{A}_{M}(Q)=\bigcup_{m\leq M}\mathcal{B}_{m}(Q),
$$
and the main result of this note can be stated as follows.

\medskip
\noindent
{\bf Theorem.} {\it Let $\pi=\otimes\pi_{v}$ and $\pi'=\otimes\pi'_v$ be in
$\mathcal{A}_{M}(Q)$. Then there exists a constant
$c=c(\varepsilon,M,K)$ depending on $\varepsilon>0, M,$ and $K$ only, such that if
$\pi_{v}\cong\pi'_{v}$ for all finite places with norm
$$
N(v)<cQ^{2M+\varepsilon}
$$
then $\pi=\pi'.$ }

\medskip

To prove the Theorem, we will exploit, among other things, Landau's classical idea
in \cite{Landau1915}.

\section{Preliminaries on automorphic $L$-functions}
\setcounter{equation}{0}
In this section, we summarize some basic properties of
automorphic $L$-functions, which will be used in the proof of the Theorem in \S3.

Let $K$ be an algebraic number
field of degree
$[K:\mathbb{Q}]=l,$
with ad\`ele ring $\Bbb A_K=K_\infty\times \Bbb A_{K,f}$, where
$K_\infty$ is the product of the Archimedean completions of $K$, and the ring $\Bbb A_{K,f}$
of finite ad\`eles is a restricted direct product of the completions $K_v$ over
non-Archimedean places $v$. Suppose that $\pi$ is an automorphic irreducible cuspidal representation of
$\GL_m(\mathbb{A}_{K}).$ Then $\pi$ is a restricted tensor product
\bea\label{pi=pipi}
\pi=\otimes_v \pi_v = \pi_\infty \otimes \pi_f,
\eea
where $v$ runs over all places of $K$, and $\pi_v$ is unramified for almost all finite
places $v$.
At every finite place $v$ where $\pi_{v}$ is unramified
we associate a semisimple
conjugacy class
\[
A_{\pi,v}=\left(
\begin{array}
[c]{ccc}%
\alpha_{\pi,v}\left(  1\right)  &  & \\
& \ddots & \\
&  & \alpha_{\pi,v}\left(  m\right)
\end{array}
\right)  ,
\]
and define the local $L$-function for the finite place $v$ as
\bea\label{eq:Euler-1}
L(s,\pi_{v})
=\det\left(I-\f{A_{\pi,v}}{q_{v}^{s}}\right)^{-1}
=\prod_{j=1}^{m}\left(1-\f{\alpha_{\pi,v}(j)}{q_{v}^{s}}\right)^{-1}
\eea
where $q_{v}=N(\mathfrak{p}_{v})=N(v)$ is the norm of $K_{v}.$
It is possible to write the local factors
at ramified places $v$ in the form
of (\ref{eq:Euler-1}) with the convention that some of the
$\alpha_{\pi,v}(j)$'s may be zero. The finite part $L$-function $L(s,\pi_f)$
is defined as
\bea
L(s,\pi_{f})=
\prod_{v<\infty}
L(s,\pi_{v}).
\eea
And this Euler product is proved to be absolutely convergent for $\s=\Re s>1$.
Also, the Archimedean $L$-function is defined as
\bea
L(s,\pi_{\infty})
=\pi^{-lms/2}
\prod_{j=1}^{lm}\Gamma\left(\frac{s+b_{\pi}(j)}{2}\right).
\eea
The coefficients $\{\alpha_{\pi,v}(j)\}_{1\leq j\leq m}$ and
$\{b_{\pi}(j)\}_{1\leq j\leq m}$ are called
local parameters of $\pi$, respectively, at finite places and at infinite places.
For them, a trivial bound states that
$$
|\alpha_{\pi,v}(j)|\leq \sqrt{p}, \qquad |\Re b_{\pi}(j)| \leq \f12.
$$
In connection with (\ref{pi=pipi}),  the complete $L$-function associated to $\pi$ is defined by
\bea
L(s,\pi)
=L(s,\pi_{\infty})
L(s,\pi_{f}).
\eea
This complete $L$-function has an analytic continuation and entire,
and satisfies the functional equation
$$
L(s,\pi)=W_{\pi}q_{\pi}^{\frac{1}{2}-s}L(1-s,\tilde\pi)
$$
where $\tilde{\pi}$ is the contragredient of $\pi,$ $W_{\pi}$ a complex
number of modulus $1$, and $q_{\pi}$ a positive integer called
the arithmetic conductor of $\pi$ \cite{Cogdell2}.

\medskip
Let $\pi=\otimes\pi_{v}$ and $\pi'=\otimes\pi'_v$
be automorphic irreducible cuspidal representations of $\mathrm{GL}_{m}(\mathbb{A}_{K})$ and
$\mathrm{GL}_{m'}(\mathbb{A}_{K})$, respectively.
The finite part Rankin-Selberg $L$-function associated to $\pi$ and $\pi'$ is defined by
\bea
L(s,\pi_{f}\times\tilde\pi'_{f})
=\prod_{v<\infty}L(s,\pi_{v}\times\tilde\pi'_{v}),
\eea
where
\bea\label{eq:coef-unramified}
L(s,\pi_{v}\times\tilde\pi'_{v})
=\prod_{j=1}^{m}\prod_{j'=1}^{m'}
\left(1-\f{\alpha_{\pi,v}(j)
\bar{\alpha}_{\pi',v}(j')}
{q_{v}^{s}}\right)^{-1}
\eea
are finite local $L$-functions for unramified finite places $v$, i.e. where
$\pi_{v}$ and $\pi'_v$ are both unramified. It can be defined similarly
at places $v$ where $\pi_{v}$ or $\pi'_v$ are ramified.
This Euler product is proved to be absolutely convergent for $\s>1$, where
$L(s,\pi_{f}\times\tilde\pi'_{f})$ has a Dirichlet series expression of the form
\bea\label{RS/DIR}
L(s,\pi_{f}\times\tilde\pi'_{f})
=\sum_{n=1}^{\infty}
\f{a_{\pi\times\tilde\pi'}(n)}{n^s}. \label{eq:coef-a}%
\eea
The complete Rankin-Selberg $L$-function is defined by
\bea
L(s,\pi\times\tilde\pi')
=L(s,\pi_{\infty}\times\tilde\pi'_{\infty})L(s,\pi_{f}\times\tilde\pi'_{f})
\eea
with
\bea\label{FORM/infty}
L(s,\pi_{\infty}\times\tilde\pi'_{\infty})
=\pi^{-mm'ls/2}\prod_{j=1}^{mm'l}
\Gamma\left(\frac{s+b_{\pi\times\tilde\pi'}(j)}{2}\right).
\eea
When the infinite place $v$ is unramified for both $\pi$ and $\pi'$, we have
$$
\{b_{\pi\times\tilde\pi'}(j)\}_{1\leq j\leq mm'}
=\{b_{\pi}(j)+b_{\tilde\pi'}(j')\}_{1\leq j\leq m, 1\leq j'\leq m'}.
$$

By Shahidi \cite{Sha1}, \cite{Sha2},
\cite{Sha3}, \cite{Sha4}, the complete $L$-function $L(s,\pi\times\tilde\pi')$
has an analytic continuation to the entire complex plane, and
satisfies the functional equation
\bea\label{FE/RSL}
L(s,\pi\times\tilde\pi')
=W_{\pi\times\tilde\pi'}
q_{\pi\times\tilde\pi'}^{\frac{1}{2}-s}
L(1-s,\tilde\pi\times\pi')
\eea
where $W_{\pi\times\tilde\pi'}$ is a complex
constant of modulus $1$, and $q_{\pi\times\tilde{\pi}^{\prime}}$ is a positive
integer. By Jacquet-Shalika \cite{JS1981} and Moeglin-Waldspurger
\cite{MoeWa1989},
we know that $L(s,\pi\times\tilde\pi')$ is holomorphic when
$\pi'\not=\pi\otimes |\det|^{i\tau}$ for any $\tau\in\mathbb{R}$, and
the only poles
of $L(s,\pi\times\tilde\pi')$ are simple poles at $s=i\tau_{0}$ and
$1+i\tau_{0}$, when $m=m'$ and
$\pi'=\pi\otimes |\det|^{i\tau_{0}}$ for some $\tau_{0}\in\mathbb{R}$.
Finally, by Gelbart-Shahidi \cite{GS2001},
$L(s,\pi\times\tilde\pi')$ is meromorphic of order one away from its poles,
and bounded in the vertical strips.

Following Iwaniec-Sarnak \cite{IS2000}, we define the analytic conductors of $L(s,\pi)$
and $L(s,\pi\times\tilde\pi')$,
respectively, as
\bea\label{ANAcon/t}
C(\pi;t)=q_{\pi}\prod_{j=1}^{ml}(1+|it+b_{\pi}(j)|),
\eea
and
\bea\label{ANAconX/t}
C(\pi,\tilde\pi';t)
=q_{\pi\times\tilde\pi'}
\prod_{j=1}^{mm'l}(1+|it+b_{\pi\times\tilde\pi'}(j)|)  .
\eea
Setting $t=0$ in the above definitions, we write
\bea\label{ANAcon}
C(\pi)=C(\pi;0), \qquad
C(\pi,\tilde\pi')
=C(\pi,\tilde\pi';0),
\eea
which are called, respectively, the analytic conductors of $\pi$ and of $\pi\times\tilde\pi'$.
They satisfy the inequality
\bea\label{eq:Bound-cond}
C(\pi,\tilde\pi')\leq C(\pi)^{m'}C(\pi')^{m},
\eea
which follows from Bushnell-Henniart \cite{BH} or Ramakrishnan-Wang \cite{RamWan}.

\section{Proof of the Theorem}
\setcounter{equation}{0}
Let $\pi=\otimes\pi_{v}$ and $\pi^{\prime}=\otimes\pi_{v}^{\prime}$
be in $\mathcal{A}_{M}(Q)$. If they are twisted equivalent, i.e.
$\pi=\pi'\otimes |\det|^{i\tau}$ for some $\tau\in {\Bbb R}^\times$,
then $\pi_v\not\cong\pi'_v$ for all finite places $v$ with at most one exception.
We may therefore suppose that,
for all $\tau\in \Bbb R$,
\bea\label{pin=pi'tau}
\pi\not=\pi'\otimes |\det|^{i\tau}. \label{eq:not-twistedequicalent}
\eea

To compare $\pi$ with $\pi'$, we form the Rankin-Selberg
$L$-function $L(s,\pi\times\tilde\pi')$, and exploit its Dirichlet series expansion
(\ref{RS/DIR}), which holds for $\s>1$. By (\ref{pin=pi'tau}), the functions
$L(s,\pi\times\tilde\pi')$ and
$L(s,\pi_f\times\tilde\pi'_f)$ are entire fucntions in the whole complex plane.
Define
\bea\label{DEF/Sxpp'}
S(x;\pi,\tilde\pi')=\sum_{n=1}^\infty a_{\pi\times\tilde\pi'}(n)w\le(\f{n}{x}\ri),
\eea
where $w(x)$ is a non-negative real valued function of $C_c^\infty$ with compact
support in $[0,3]$, and we may specify
$$
w(x)=
\left\{
\begin{array}{lll}
0, & \text{for } x\not\in [0,3], \\
e^{-\f1x}, & \text{for } x\in (0,1], \\
e^{-\f{1}{3-x}}, & \text{for } x\in [2,3)
\end{array}
\right.
$$
as in \cite[\S3]{Wang}. Thus, for any positive integer $k$, the derivative
$w^{(k)}(x)$ has exponential decay as $x\to 0$ or $3$. Consequently, the Mellin transform
$$
W(s)=\int_0^\infty w(x)x^{s-1}dx
$$
is an analytic function of $s$; if $\s<-1$ then
$$
W(s)\ll _{A,\s}\f{1}{(1+|t|)^A}
$$
for arbitrary $A>0$, by repeated partial integration. By Mellin inversion,
$$
w(x)=\f{1}{2\pi i}\int_{(2)} W(s)x^{-s}ds,
$$
where $(c)$ means the vertical line $\s=c$.
Inserting this back to (\ref{DEF/Sxpp'}), and then using Dirichlet series expansion
(\ref{RS/DIR}), we have
\bna\label{DEF/Sxpp'+}
S(x;\pi,\tilde\pi')
&=&\f{1}{2\pi i}\sum_{n=1}^\infty a_{\pi\times\tilde\pi'}(n)
\int_{(2)} W(s)\le(\f{n}{x}\ri)^{-s}ds\nonumber\\
&=&\f{1}{2\pi i}\int_{(2)} x^s W(s) L(s,\pi_f\times\tilde\pi'_f)ds,
\ena
where the interchange of summation and integral is guaranteed by the absolute convergence
of (\ref{RS/DIR}) on the line $\s=2$. A pre-convexity bound like
$$
L(s,\pi_f\times\tilde\pi'_f)\ll C(\pi,\tilde\pi';t)^B,
$$
where $B>0$ is some constant, can be obtained by standard method, as pointed out
in \cite[\S1]{Brumley2004}. Since $W(s)$ is rapid decay and
$L(s,\pi\times\tilde\pi')$ is entire, we may
shift the contour above to the vertical line $\s=-H$, getting
\bea\label{DEF/Sxpp'+}
S(x;\pi,\tilde\pi')
=\f{1}{2\pi i}\int_{(-H)}x^sW(s)L(s,\pi_f\times\tilde\pi'_f)ds,
\eea
where $H>1$ is a large real number to be specified later.

We are going to apply the functional equation (\ref{FE/RSL}) to (\ref{DEF/Sxpp'+}). To this
end, we rewrite (\ref{FE/RSL}) as
\bea\label{FE/REWRITE}
L(s,\pi_f\times\tilde\pi'_f)=W_{\pi\times\tilde\pi'}q_{\pi\times\tilde\pi'}^{\f12-s}
G(s)L(1-s,\tilde\pi_f\times\pi'_f),
\eea
where
\bea\label{def/G}
G(s)=\f{L(1-s,\tilde\pi_\infty\times\pi'_\infty)}{L(s,\pi_\infty\times\tilde\pi'_\infty)}.
\eea
We need to estimate $G(s)$ on the line $\s=-H$, avoiding the poles of
the nominator of $G(s)$. This will be done in the Lemma in \S4,
which asserts that, for every large positive integer $n$, there is an $H\in [n,n+1]$
such that, on the vertical line $\s=-H$,
\bea\label{G(s)<<}
G(-H+it)
\ll_{H,M,K} (1+|t|)^{mm'l(\f12+H)}\prod_{j=1}^{mm'l}(1+|b_{\pi\times\tilde\pi'}(j)|)^{\f12+H}.
\eea
Now we apply the functional equation (\ref{FE/REWRITE}),
\bna
S(x;\pi,\pi')
&=&\f{1}{2\pi i}\int_{(-H)}x^sW(s)W_{\pi\times\tilde\pi'}q^{\f12-s}_{\pi\times\tilde\pi'}
G(s)L(1-s,\tilde\pi_{f}\times \pi'_{f})ds \nonumber\\
&=&\f{1}{2\pi i}\int_{(-H)}x^sW(s)W_{\pi\times\tilde\pi'}q^{\f12-s}_{\pi\times\tilde\pi'}
G(s)\le(\sum_{n=1}^\infty\f{a_{\tilde\pi\times\pi'}(n)}{n^{1-s}}\ri)ds \nonumber\\
&=&\f{1}{2\pi i}\sum_{n=1}^\infty\f{a_{\tilde\pi\times\pi'}(n)}{n^{1+H}}
\int_{(-H)}x^sW(s)W_{\pi\times\tilde\pi'}q^{\f12-s}_{\pi\times\tilde\pi'}
G(s)n^{s+H}ds.
\ena
Here the interchange of summation and integral is guaranteed by
the absolute convergence of the Dirichlet series as well as the rapid
decay of $W(s)$. Using these facts again, and inserting (\ref{G(s)<<}) into the last integral, we get
\bea\label{Sxpp'=}
S(x;\pi,\tilde\pi')
&\ll_{H,M,K}& \int_{(-H)}\bigg|x^sW(s)W_{\pi\times\tilde\pi'}q^{\f12-s}_{\pi\times\tilde\pi'}
G(s)ds\bigg|\nonumber\\
&\ll_{H,M,K}& x^{-H} q^{\f12+H}_{\pi\times\tilde\pi'} \prod_{j=1}^{mm'l}(1+|b_{\pi\times\tilde\pi'}(j)|)^{\f12+H}
\nonumber\\
&=& x^{-H}C(\pi,\tilde\pi')^{\f12+H}.
\eea
This upper bound corresponds to \cite[(10)]{Wang}.

To establish the Theorem, we need a lower bound for $S(x;\pi,\tilde\pi')$.
Under the assumption of (\ref{eq:not-twistedequicalent}), we further suppose that $\pi_{v}\cong\pi'_{v}$ for all finite places with norm
$N(v)<x$. Then
\bea\label{S=S}
S(x;\pi,\tilde\pi')
=S(x;\pi,\tilde\pi).
\eea
A lower bound for $S(x;\pi,\tilde\pi)$ is obtained in
\cite[Lemma 2]{Brumley2004}; see also \cite[Propsition 4]{Wang}. Thus,
we have
\bea\label{S>1}
S(x;\pi,\tilde\pi)\gg \f{x^{\f1m}}{\log x}.
\eea

Combining (\ref{Sxpp'=}), (\ref{S=S}), and (\ref{S>1}), we get
$$
\f{x^{\f1m}}{\log x} \ll S(x;\pi,\tilde\pi')\ll_{H,M,K} x^{-H}C(\pi,\tilde\pi')^{\f12+H}.
$$
One therefore has
$$
x<c\cdot C(\pi,\tilde\pi')^{\f{H+1/2}{H+1/m}},
$$
where $c$ is a constant depending on $H,M,K$. Taking $H$ sufficiently large,
this becomes
$$
x<c\cdot C(\pi,\tilde\pi')^{1+\ve},
$$
and now the constant $c$ depends on $\ve,M,K$. The assertion of the theorem finally
follows from this and (\ref{eq:Bound-cond}).

\section{An Estimate for $G(s)$}
\setcounter{equation}{0}

In this section, we give an estimate for $G(s)$ defined as
in (\ref{def/G}) on a vertical line $s=-H+it$, where
$H$ is a large real number to be decided suitably.
Recall Stirling's formula that
$$
|\G(\s+it)|=\sqrt{2\pi} e^{-\f{\pi}{2}|t|}|t|^{\s-\f12}
\bigg(1+O_{\s,\d}\bigg(\f{1}{1+|t|}\bigg)\bigg),
$$
which holds for $s=\s+it$ away from all poles of $\G(s)$ by at least $\d>0$;
note that the implied constant depends on $\s$ and $\d$.

To do this, we should first locate the poles of the nominator of $G(s)$, i.e.
poles of
\bea\label{L1-s=G}
L(1-s,\tilde\pi_\infty\times\pi'_\infty)
=\pi^{-mm'ls/2}\prod_{j=1}^{mm'l}
\Gamma\left(\frac{1-s+b_{\tilde\pi\times\pi'}(j)}{2}\right),
\eea
according to (\ref{FORM/infty}). These poles are easily to be seen as
$$
P_{n,j}=2n+1+b_{\tilde\pi\times\pi'}(j), \quad n=0,1,2,\cdots, \quad j=1,\cdots,mm'l.
$$
As in \cite{LiuYe}, we let ${\Bbb C}(m,m')$ be the complex plane with the discs
$$
|s-P_{n,j}|<\f{1}{8mm'l},
\quad n=0,1,2,\cdots, \quad j=1,\cdots,mm'l
$$
excluded. Thus, for any $s\in {\Bbb C}(m,m')$, the quantity
$$
\frac{1-s+b_{\tilde\pi\times\pi'}(j)}{2}
$$
is away from all poles of  $\Gamma(s)$ by at least $1/(16mm'l)$, and therefore
Stirling's formula applies to (\ref{L1-s=G}). Of course, Stirling's
formula also applies to the denominator of $G(s)$.  Writing
$b_{\pi\times\tilde\pi'}(j)=u(j)+iv(j)$ and $s=\s+it\in \Bbb C(m,m')$, we have
\bna
G(s)
&=&\pi^{-\f{mm'l}{2}+mm'ls}\prod_{j=1}^{mm'l}
\f{\Gamma\left(\frac{1-s+b_{\tilde\pi\times\pi'}(j)}{2}\right)}
{\Gamma\left(\frac{s+b_{\pi\times\tilde\pi'}(j)}{2}\right)} \\
&\ll_{\s,M,K}&
\prod_{j=1}^{mm'l} \f{|t+v(j)|^{\f{1-\s+u(j)}{2}-\f12}}
{|t+v(j)|^{\f{\s+u(j)}{2}-\f12}} \\
&\ll_{\s,M,K}&
\prod_{j=1}^{mm'l} |t+v(j)|^{\f12-\s},
\ena
where we have used $\{b_{\tilde\pi\times\pi'}(j)\}_{1\leq j\leq mm'l}
=\{\overline{b_{\pi\times\tilde\pi'}(j)}\}_{1\leq j\leq mm'l}$.
It follows that, for $\s<1/2$,
\bea
G(s)\ll_{\s,M,K}
(1+|t|)^{mm'l(\f12-\s)} \prod_{j=1}^{mm'l} (1+|v(j)|)^{\f12-\s},
\eea
which can be written as
\bea\label{EST/G(s)}
G(s)\ll_{\s,M,K}
(1+|t|)^{mm'l(\f12-\s)} \prod_{j=1}^{mm'l} (1+|b_{\pi\times\tilde\pi'}(j)|)^{\f12-\s}.
\eea

Now we give a remark about the structure of ${\Bbb C}(m,m')$. For
$j=1,\cdots,mm'l$, denote by $\beta(j)$ the fractional
part of $v(j)$. In addition we let
$\beta(0)=0$ and $\beta(mm'l+1)=1$.
Then all $\beta(j)\in [0,1]$, and hence there exist
$\beta(j_1),\beta(j_2)$ such that
$\beta(j_2)-\beta(j_1)\geq 1/(3mm'l)$ and there is no
$\beta(j)$ lying between $\beta(j_1)$ and $\beta(j_2)$. It
follows that the strip
$S_0=\{s:\beta(j_1)+1/(8mm'l)\leq \Re s\leq\beta(j_2)-1/(8mm'l)\}$
is contained in ${\Bbb C}(m,m').$ Consequently, for all $n=0,1,2,\cdots,$
the strips
$$
S_n=\bigg\{s: -n+\beta(j_1)+\f{1}{8mm'l}\leq \Re s\leq -n+\beta(j_2)-\f{1}{8mm'l}\bigg\}
$$
are subsets of ${\Bbb C}(m,m').$ Therefore, for each $n\geq 1$, one can
choose a vertical line $\s=-H$ lying in the strip $S_n$, and therefore
(\ref{EST/G(s)}) holds on the line $\s=-H$. This proves the following result.

\medskip
\noindent
{\bf Lemma.} {\it Let $G(s)$ be as in (\ref{def/G}). Then for each $n\geq 1$, there
is an $H\in [n,n+1]$, such that on the line $\s=-H$ we have
\bna
G(-H+it)\ll_{H,M,K}
(1+|t|)^{mm'l(\f12+H)} \prod_{j=1}^{mm'l} (1+|b_{\pi\times\tilde\pi'}(j)|)^{\f12+H}.
\ena}

\bigskip
\bigskip

\footnotesize \noindent \halign{#\hfill \quad           & #\hfill
\quad & #\hfill \cr Jianya Liu                      & Yonghui Wang
\cr School of Mathematics           & Department of Mathematics \cr
Shandong University             & Capital Normal University \cr
Jinan, Shandong 250100          & Beijing 100037 \cr China & China
\cr {\tt jyliu@sdu.edu.cn}          & {\tt arith.yonghui.wang@gmail.com} \cr
                                & {\tt yhwang@mail.cnu.edu.cn}                \cr}
\end{document}